\input amstex
\input epsf
\magnification=\magstep1 
\baselineskip=13pt
\documentstyle{amsppt}
\vsize=8.7truein \CenteredTagsOnSplits \NoRunningHeads
\def\ttt{\bold{t}}
\def\xx{\bold{x}}
\def\EE{\bold {E\thinspace}}
\def\tr{\operatorname{trace}}
\def\rk{\operatorname{rank}}
\def\dist{\operatorname{dist}}
\def\supp{\operatorname{supp}}

\def\GG{\Cal{G}}

\topmatter
 
\title Testing systems of real quadratic equations for approximate solutions   \endtitle 
\author Alexander Barvinok  \endauthor
\address Department of Mathematics, University of Michigan, Ann Arbor,
MI 48109-1043, USA \endaddress
\email barvinok$\@$umich.edu \endemail
\date June 23,  2020 \enddate
\thanks  This research was partially supported by NSF Grant DMS 1855428. 
\endthanks 
\abstract Consider systems of equations $q_i(x)=0$, where $q_i: {\Bbb R}^n \longrightarrow {\Bbb R}$, $i=1, \ldots, m$, are quadratic forms. Our goal is to tell efficiently systems with many non-trivial solutions or near-solutions $x \ne 0$ from systems that are far from having a solution. For that, we pick a delta-shaped penalty function $F: {\Bbb R} \longrightarrow [0, 1]$ with $F(0)=1$ and $F(y) < 1$ for $y \ne 0$
and compute the expectation of $F(q_1(x)) \cdots F(q_m(x))$ for a random $x$ sampled from the standard Gaussian measure in ${\Bbb R}^n$. We choose $F(y)=y^{-2}\sin^2 y$ and show that the expectation can be approximated within relative error $0< \epsilon < 1$ in quasi-polynomial time $(m+n)^{O(\ln (m+n)-\ln \epsilon)}$, provided each form $q_i$ depends on not more than $r$ real variables, has common variables with at most $r-1$ other forms and satisfies $|q_i(x)| \leq \gamma \|x\|^2/r$, where $\gamma >0$ is an absolute constant. This allows us to distinguish between ``easily solvable" and ``badly unsolvable" systems in some non-trivial situations.
 \endabstract
\endtopmatter
\document

\head 1. Introduction and main results \endhead

\subhead (1.1) Systems of real quadratic equations \endsubhead 
We consider systems of homogeneous real quadratic equations
$$q_i(x)=0\quad \text{for} \quad  i=1, \ldots, m, \tag1.1.1$$
where $q_i: {\Bbb R}^n \longrightarrow {\Bbb R}$ are quadratic forms,
$$q_i(x)={1 \over 2} \langle Q_i x, x \rangle \quad \text{for} \quad i=1, \ldots, m. $$
Here $\langle \cdot, \cdot \rangle$ is the standard scalar product in ${\Bbb R}^n$ and $Q_i$ are self-adjoint operators, represented by $n \times n$ symmetric matrices in the standard
basis of ${\Bbb R}^n$. We are interested in finding out whether the system (1.1.1) has a non-trivial solution $x \ne 0$. Generally, the problem is computationally hard (it is not even known to be in NP, since the description of a solution $x$ can have exponential complexity). However, if the number 
$m$ of equations is fixed in advance, a polynomial time algorithm is available \cite{Ba93}, see also \cite{GP05} and \cite{Ba08} for strengthening (but as a function of $m$, the complexity of the algorithm is exponential). If the number $n$ of variables is fixed in advance, the feasibility of the system and many related problems (also when $q_i$ are arbitrary polynomials, not necessarily quadratic) can be solved in polynomial time, but as a function of $n$, the complexity of the algorithm is exponential, see for example, 
\cite{B+06}.

We note that finding whether any given system of real polynomial equations can be reduced to finding whether a system of homogeneous quadratic equations has a non-trivial solution. First, by the introduction of new variables and repeated substitutions of the type $xy =z$, we successively lower the degree of polynomials. This way we arrive to a system of equations of the type $q_i(x) =0$, where $q_i$ are quadratic polynomials. Then we introduce yet another variable $t$ and replace each equation $q_i(x)=0$ by $t^2 q_i\left(t^{-1} x \right)=0$, making all equations homogeneous. It remains to make sure that $t \ne 0$ by introducing an additional variable $s$ and the equation $R^2 t^2 - \left(x_1^2 + \ldots + x_n^2\right) = s^2$ binding all variables together, so that if $t=0$ then all other variables are also $0$ (here $R$ is treated as a very large constant, in fact it can be treated as infinitely large, with computations in the field of rational functions in $R$, cf. \cite{GV88} ).

We also note that in some areas, for example in distance geometry, systems of quadratic equations appear naturally, see \cite{L+14}. In particular, the equations that appear in distance geometry tend to be sparse: there the unknowns are $d$-vectors $x_1, \ldots, x_n \in {\Bbb R}^d$ and the equations relate squared Euclidean distances $\|x_i - x_j \|^2$ between pairs of points, each of which is a quadratic form in $2d$ real variables.
Furthermore, the number of equations involving a particular variable is determined by the number of conditions imposed on a particular vector.

\subhead (1.2) Accounting for solutions \endsubhead Since testing the feasibility of (1.1.1) in the general case is computationally hard, we pursue a more modest goal. We would like to be able to efficiently separate 
the systems that have ``many near-solutions" from the systems that are ``far from having a solution".
To accomplish our goal, we introduce a ``penalty function" $F: {\Bbb R} \longrightarrow [0, 1]$ such that 
$F(0)=1$ and $F(y) < 1$ for $y \ne 0$,  and compute the integral 
$${1 \over (2 \pi)^{n/2}} \int_{{\Bbb R}^n} F\left(q_1(x)\right) \cdots F\left(q_m(x)\right) e^{-\|x\|^2/2} \ dx. \tag1.2.1$$
The goal is to choose $F$ as sharply peaked at 0 as possible, so that a point $x \in {\Bbb R}^n$ is accounted for with the Gaussian weight $(2 \pi)^{-n/2} e^{-\|x\|^2/2}$, if $x$ is a solution of (1.1.1) and is accounted for with an exponentially smaller weight if for many of the forms $q_i$, the values of $q_i(x)$ are far from $0$. Hence we expect the integral (1.2.1) to be large if there are many solutions or ``near-solutions" $x$ and we expect (1.2.1) to be small if the system is far from having a solution, in which case even a small perturbation $q_i \longmapsto \tilde{q}_i$ does not result in a system having a non-trivial solutions. 

We are able to choose 
$$F(y)={\sin^2 y \over y^2},$$
where $F(0)=1$ by continuity.
\epsfxsize=2truein 
$$\epsffile{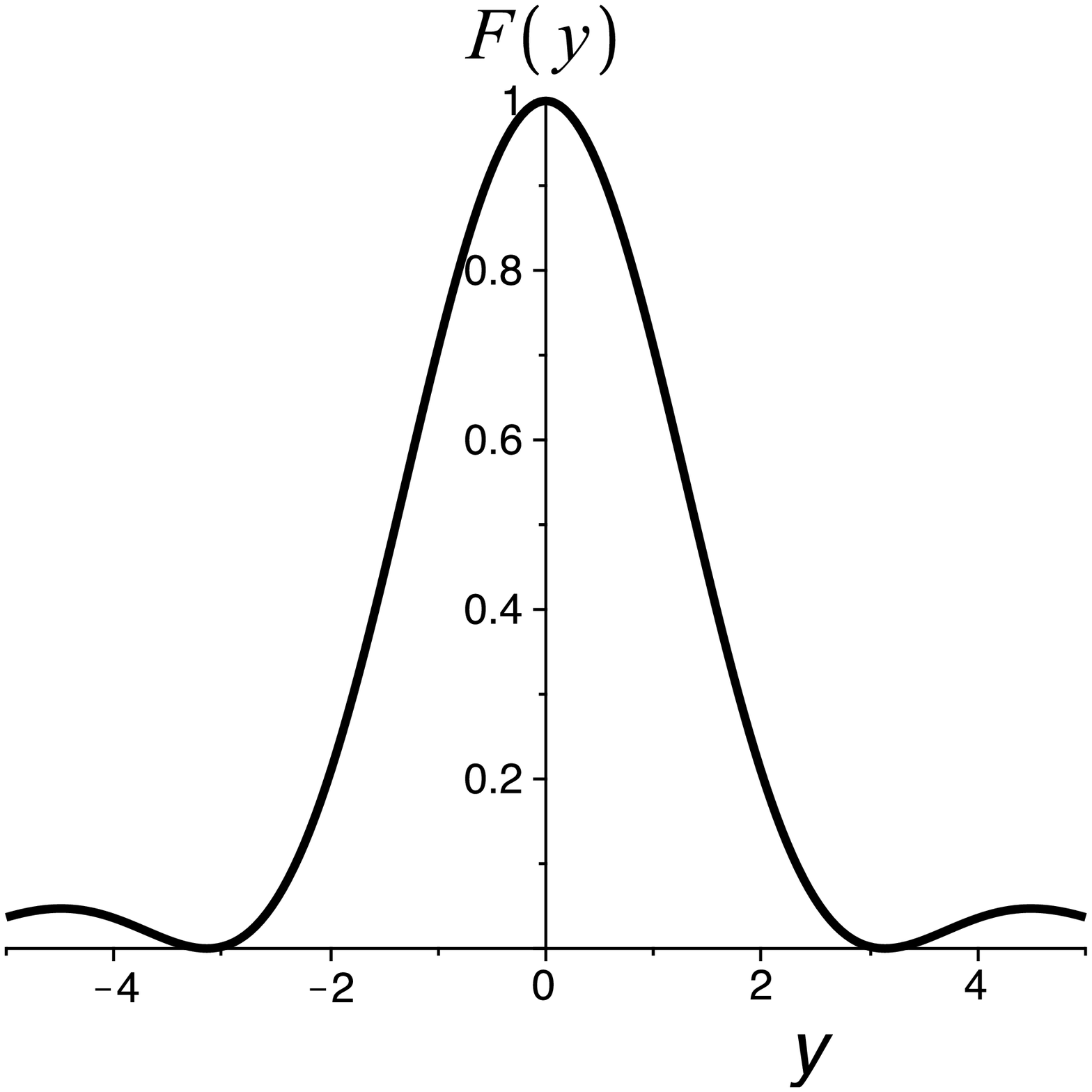}$$
Let
$$\|x\|=\sqrt{\langle x, x \rangle} =\sqrt{x_1^2 + \ldots + x_n^2} \quad \text{for} \quad x=\left(x_1, \ldots, x_n\right)$$
denote the standard Euclidean norm in ${\Bbb R}^n$. We prove the following main result.
\proclaim{(1.3) Theorem} There is an absolute constant $\gamma_1 >0$ (one can choose $\gamma_1=0.09$) such that the following holds. 
Let $q_i: {\Bbb R}^n \longrightarrow {\Bbb R}$, $i=1, \ldots, m$, be quadratic forms in $n$ real variables $x_1, \ldots, x_n$, such that each form $q_i$ depends on at most $r$ variables 
among $x_1, \ldots, x_n$, has common variables with at most $r-1$ other forms $q_j$ and satisfies
$$|q_i(x)| \ \leq \ {\gamma_1 \|x\|^2 \over r} \quad \text{for} \quad i=1, \ldots, m.$$
Then for any $0 < \epsilon < 1$, one can compute the value of 
$${1 \over (2 \pi)^{n/2}} \int_{{\Bbb R}^n} e^{-\|x\|^2/2} \prod_{i=1}^m { \sin^2 \left(q_i(x)\right) \over q_i^2(x)} \ dx \tag1.3.1$$
within relative error $\epsilon$ in quasi-polynomial $(m+n)^{O(\ln (m+n) - \ln \epsilon)}$ time.
\endproclaim

 As is well-known, the bulk of the standard Gaussian measure $\mu$ in ${\Bbb R}^n$ with density $(2\pi)^{-n/2} e^{-\|x\|^2/2}$ is concentrated in the vicinity of the sphere $\|x\|=\sqrt{n}$. More precisely, 
$$\mu\left\{x: \ (1-\epsilon) n \ \leq \  \|x\|^2 \ \leq \ {n \over 1-\epsilon} \right\} \ \geq \ 1- 2 e^{-\epsilon^2 n/4} \quad \text{for} \quad 0 < \epsilon < 1,$$
see, for example, Section V.5 of \cite{Ba02}. Under the conditions of Theorem 1.3, for a random point $x \in {\Bbb R}^n$, we can expect the values of $q_i(x)$ to be of the order of a non-zero constant.

Let us consider an asymptotic regime where $m=o(n)$ and both $m$ and $n$ grow, while $r$ may grow also.
Hence if the value of (1.3.1) is $\beta^m$ for some $0 < \beta < 1$, we are guaranteed to have a point $x_0 \in {\Bbb R}^n$ with $\|x_0\| = \sqrt{n}(1+o(1))$ and 
$$\prod_{i=1}^m {\sin^2 \bigl(q_i(x_0)\bigr) \over q_i^2(x_0)} \geq \beta^m.$$
It is not clear though how to construct such a point $x_0$ efficiently. Another question is whether for $r$, fixed in advance, one can approximate (1.3.1) in genuine polynomial, as opposed to quasi-polynomial time, 
building perhaps on the methods of \cite{PR17}.

Below, we discuss some situations where the algorithm of Theorem 1.3 allows us to tell, in some asymptotic regimes, systems with ``many solutions", where the value of (1.3.1) is large, from systems that are ``far from having a solution",
where the value of (1.3.1) is small.

\subhead (1.4) Separating systems with many solutions from systems that are far from having a solution \endsubhead 
Suppose that in (1.1.1) we have $n=k m$ for some integer $k$. First, we consider the case of general quadratic forms $q_i: {\Bbb R}^n \longrightarrow {\Bbb R}$, $i=1, \ldots, m$. Scaling, if necessary, we assume that 
$$|q_i(x)| \ \leq \ {\gamma_1 \|x\|^2 \over n}, \tag1.4.1$$
where $\gamma_1$ is a constant in Theorem 1.3. Next, we introduce $k$ copies of each quadratic form $q_i$ and consider the integral (1.3.1), which in our case is written as 
$${1 \over (2 \pi)^{n/2}} \int_{{\Bbb R}^n} e^{-\|x\|^2/2} \prod_{i=1}^m \left({\sin^2 \bigl(q_i(x) \bigr)\over q_i^2(x)}\right)^k \ dx. \tag1.4.2$$
Using the algorithm of Theorem 1.3, we can approximate (1.4.2) within relative error $0 < \epsilon < 1$ in quasi-polynomial time $n^{O(\ln n - \ln \epsilon)}$.

Suppose first, that there are many solutions, by which we mean the following.
Let 
$$X=\Bigl\{x \in {\Bbb R}^n: \quad q_i(x)=0 \quad i=1, \ldots, m \Bigr\}$$
be the set of solutions. In general position, if non-empty, $X$ is a real algebraic variety of codimension $m$. 
For $\delta >0$, let
$$X_{\delta}=\Bigl\{y \in {\Bbb R}^n:\quad \dist(y, X) \leq \delta \Bigr\}$$
be the $\delta$-neighborhood of $X$, where 
$$\dist(y, X)=\min_{x \in X} \|y-x\|$$
is the distance from $y$ to $X$. We say that (1.1.1) has ``many solutions" if for the standard Gaussian measure $\mu$ and for some fixed $\kappa \geq 2$, we have 
$$\mu\left(X_{n^{-\kappa}}\right) = n^{-O(m)}  \tag1.4.3$$
(we assume that all implied constants in the ``$O$" notation are absolute). This is the case, for example, if $X$ is a subspace of codimension $m$, or if $n=2d$ is even, ${\Bbb R}^n$ is identified with ${\Bbb C}^d$ and 
$q_i$ are identified with complex quadratic polynomials, see \cite{Kl18}, \cite{AK18}. It follows from the Gromov theorem on the waist of the Gaussian measure \cite{Gr03} that for any set of quadratic forms as in (1.1.1)
one can find $\alpha_1, \ldots, \alpha_m$ such that (1.4.3) holds for the modified forms $\widehat{q}_i(x)=q_i(x)-\alpha_i \|x\|^2$.

We consider the asymptotic regime in which $n \longrightarrow \infty$ and $k \gg \ln m$, so that $m=o(n/\ln n)$ (and where $m$ may also grow). It is not hard to see that if (1.4.3) holds, the value of (1.4.2) is $n^{-O(m)}$.   
This follows since the bulk of the measure $\mu$ is concentrated in the vicinity of the sphere $\|x\|=\sqrt{n}$, where the Lipschitz constant of each quadratic form $q_i$ is $O(1/\sqrt{n})$ and hence the value of the product 
$$\prod_{i=1}^m \left({\sin^2 \bigl(q_i(x) \bigr)\over q_i^2(x)}\right)^k$$ in 
(1.4.2) for each $x \in X_{n^{-\kappa}}$ is $1-o(1)$.

Suppose now that under the same asymptotic regime, the system (1.1.1) is far from having a solution, by which we mean that for some fixed $\delta >0, \beta >0$ and all $x \in {\Bbb R}^n$ such that $\|x\|=\sqrt{n}$, for at least $\delta m$ of the quadratic forms $q_i$, we have $|q_i(x)| \geq \beta$. In that case, we conclude that the value of (1.4.2) is $2^{-\Omega(n)}$ and we can tell apart the two cases as $n$ grows.

One can argue that in this particular situation there is a simple randomized algorithm telling the two cases apart: one should pick a random point $x \in {\Bbb R}^n$ with respect to the standard Gaussian distribution. If the system has many solutions (in the above sense), then with high probability one should have 
$$\max_{i=1, \ldots, m} |q_i(x)| = o\left(1\right)$$
and if the system is far from having a solution (again, in the above sense) then this condition is violated. This follows from the isoperimetric inequality for the Gaussian measure, see, for example, Section 4.3 of \cite{Bo98} and the fact that in the vicinity of the sphere $\|x\|=\sqrt{n}$, where the Gaussian measure is concentrated, the Lipschitz constants of the forms $q_i$ are all $O(n^{-1/2})$.
However, there appears to be no equally obvious deterministic algorithm that would tell the difference. Moreover, as we venture into sparse systems, where each form $q_i$ depends on at most $rk$ variables and has a common variable with at most $r-1$ other forms , we replace the condition (1.4.1) by a weaker condition
$$|q_i(x)| \ \leq \ {\gamma_1 \|x\|^2 \over rk} \quad \text{for} \quad i=1, \ldots, m$$
and although our approach still works, no other randomized or deterministic algorithm seems to be able to tell efficiently the two cases apart.

\head 2. The plan of the proof \endhead

The algorithm of Theorem 1.3 is based on the method of polynomial interpolation, see Section 2.2 of \cite{Ba16} and also \cite{PR17} for some enhancement. The idea is to construct a univariate polynomial $\phi_N(z)$ of some degree 
$N=O\left((m+n)\ln (1/\epsilon)\right)$ such that the following conditions are satisfied:
\bigskip
$\bullet$ The value of $\phi_N(1)$ approximates the value of (1.3.1) within the desired relative error $\epsilon$;
\smallskip
$\bullet$ For some $\beta >1$, we have $\phi_N(z) \ne 0$ for all $z \in {\Bbb C}$ satisfying $|z| < \beta$ and
\smallskip
$\bullet$ For any $k$, one can compute the derivative $\phi_N^{(k)}(0)$ in $N^{O(k)}$ time.
\bigskip
As discussed in Section 2.2 of \cite{Ba16}, one can then approximate $\phi_N(1)$ within relative error $0 < \epsilon < 1$ from the values of $\phi_N^{(k)}(0)$ for $k=O(\ln N - \ln \epsilon)$ and hence one can approximate 
$\phi_N(1)$ in $N^{O(\ln N -\ln \epsilon)}$ time.

In our case, we define $\phi_N(z)$ as the Taylor polynomial of a sufficiently large degree $N$ of the function 
$$\phi(z)={1 \over (2\pi)^{n/2}} \int_{{\Bbb R}^n} e^{-\|x\|^2/2} \prod_{i=1}^m {\sin^2 \bigl(z q_i(x) \bigr) \over z^2 q_i^2(x)} \ dx,$$
so that $\phi_N(1) \approx \phi(1)$. To ascertain that the polynomial $\phi_N(z)$ satisfies the required properties, we prove the following result.

\proclaim{(2.1) Theorem} There is an absolute constant $0 < \gamma < 0.25$ (one can choose $\gamma=0.1$) such that the following holds. Let $q_i: {\Bbb R}^n \longrightarrow {\Bbb R}$, $i=1, \ldots, m$,  
be quadratic forms in real variables $x_1, \ldots, x_n$, such that each form $q_i$ depends on at most $r$ variables among $x_1, \ldots, x_n$, has common variables with at most $r-1$ other forms $q_j$ and 
satisfies 
$$|q_i(x)| \ \leq \ {\gamma \|x\|^2 \over r} \quad \text{for} \quad i=1, \ldots, m.$$ Then, for all $z \in {\Bbb C}$ such that $|z| \leq 1$, the integral 
$$\phi(z)={1 \over (2 \pi)^{n/2}} \int_{{\Bbb R}^n} e^{-\|x\|^2/2} \prod_{i=1}^m {\sin^2 \left(z q_i(x)\right) \over z^2 q_i^2(x)} \ dx$$
converges absolutely to an analytic function. Moreover,
$$(1-4 \gamma)^{-n/2} \ \geq \ |\phi(z)| \ \geq \ 2^{-m/2} (1+4 \gamma)^{-n/2} \quad \text{for all} \quad |z| \leq 1.$$
\endproclaim

We prove Theorem 2.1 in two steps. First, we find a different integral representation for $\phi(z)$.

Let 
$$[-1, 1]^m=\Bigl\{\left(t_1, \ldots, t_m\right): \ -1 \leq t_i \leq 1 \quad \text{for} \quad i=1, \ldots, m \Bigr\}$$ 
be the cube endowed with Lebesgue measure $dt$. For an $n \times n$ complex $Q$, we denote by $\|Q\|$ its operator norm, that is, 
$$\|Q\|=\max_{x \in {\Bbb C}^n: \ \|x\|=1} \|Q x\|.$$ We denote by $\sqrt{-1}$ the imaginary unit, so as to save ``$i$" for indexing.

\proclaim{(2.2) Proposition} Let $q_1, \ldots, q_m: {\Bbb R}^n \longrightarrow {\Bbb R}$ be quadratic forms,
$$q_i(x)={1 \over 2} \langle Q_i x, x \rangle \quad \text{for} \quad i=1, \ldots, m,$$
where $Q_1, \ldots, Q_m$ are $n \times n$ real symmetric matrices satisfying
$$ \left\| \sum_{i=1}^m z_i Q_i \right\| \ < \ {1 \over 2}$$
for all $z_1, \ldots, z_m \in {\Bbb C}$ satisfying $|z_1|, \ldots, |z_m| \leq 1$.
Then for all $z \in {\Bbb C}$ such that $|z| \leq 1$, we have 
$$\split &{1 \over (2 \pi)^{n/2}} \int_{{\Bbb R}^n} e^{-\|x\|^2/2} \prod_{i=1}^m {\sin^2 \bigl(z q_i(x)\bigr) \over z^2 q_i^2(x)} \ dx \\ &\qquad =\int_{[-1, 1]^m} \prod_{i=1}^m \left(1-|t_i|\right)
 \det^{\qquad -{1 \over 2}} \left( I - 2z\sum_{i=1}^m t_i \sqrt{-1} Q_i \right) \ dt. \endsplit$$
\endproclaim 

Next, we expand 
$$\det^{\qquad -{1 \over 2}} \left( I - 2z\sum_{i=1}^m   t_i \sqrt{-1} Q_i \right) = \exp\left\{ \sum_{a \in {\Bbb Z}^m_+} c_a \ttt^a  \right\},$$
where ${\Bbb Z}^m_+$ is the set of non-negative integer $m$-vectors and
$$\ttt^a=t_1^{\alpha_1} \cdots t_m^{\alpha_m} \quad \text{for} \quad a=\left(\alpha_1, \ldots, \alpha_m\right),$$
and where we have $c_a = c_a(z) \in {\Bbb C}$. 

For a vector $a \in {\Bbb Z}^m_+$, $a=\left(\alpha_1, \ldots, \alpha_m\right)$, we define the {\it support} of $a$ as 
$$\supp a =\left\{ i:\ \alpha_i \ne 0 \right\}.$$
To complete the proof, we obtain the following fairly general result (in fact, we prove it in more generality than we need in this paper).
\proclaim{(2.3) Theorem} Let $\nu$ be a product measure on the cube $[-1, 1]^m$, so $\nu=\nu_1 \times \ldots \times \nu_m$, where $\nu_i$ is a measure on the $i$-th copy of the interval $[-1, 1]$.
Let 
$$g(t_1, \ldots, t_m)=\sum_{a \in {\Bbb Z}^m_+} c_{a} \ttt^{a}, $$
and $c_a \in {\Bbb C}$ for all $a \in  {\Bbb Z}^m_+$ and the series converges absolutely and uniformly on $[-1, 1]^m$.
Suppose that for some $0 \leq \theta_1, \ldots, \theta_m < 2\pi/3$ and all $i=1, \ldots, m$, we have 
$$\sum\Sb a \in {\Bbb Z}^m_+: \\ i \in \supp a \endSb |c_a| \prod_{j \in \supp a} {1 \over \cos (\theta_j/2)} \ \leq \ {\theta_i \over 2}.$$
Then 
$$\left| \int_{[-1, 1]^m} e^g \ d \nu \right| \ \geq \ \left(\prod_{i=1}^m  \cos(\theta_i/2) \right) \int_{[-1, 1]^m} \left| e^{g}\right| \ d \nu.$$
\endproclaim

We prove Proposition 2.2 in Section 3, Theorem 2.3 in Section 4 and Theorem 2.1 in Section 5. In Section 6, we prove Theorem 1.3 by providing the algorithm to approximate (1.3.1).

\head 3. Proof of Proposition 2.2 \endhead

We start with some simple calculations.
\proclaim{(3.1) Lemma} We have 
$$2 \int_0^1 (1-t) \cos (2yt) \ dt = {\sin^2 y \over y^2}.$$
\endproclaim
\demo{Proof} The formula obviously holds for $y=0$. For $y \ne 0$, integrating by parts, we obtain
$$\split &2 \int_0^1 (1-t) \cos (2yt) \ dt= y^{-1} \int_0^1 (1-t) \ d \left(\sin (2yt) \right) \\ &\qquad= (1-t){\sin (2yt) \over y} \Big|_{t=0}^{t=1}  + y^{-1} \int_0^1 \sin (2 yt ) \ dt =
- {\cos (2y t) \over 2y^2} \Big|_{t=0}^{t=1} \\ 
&\qquad={1 -\cos (2y) \over 2y^2}={\sin^2 y \over y^2}. \endsplit $$
{\hfill \hfill \hfill} \qed
\enddemo

\proclaim{(3.2) Lemma} Let $q: {\Bbb R}^n \longrightarrow {\Bbb R}$ be a quadratic form, 
$$q(x)= {1 \over 2} \langle Qx, x \rangle,$$
where $Q$ is an $n \times n$ symmetric matrix such that $\| Q\| \leq 1$. Then for all $z \in {\Bbb C}$ such that $|z| <1$, we have 
$${1 \over (2\pi)^{n/2}} \int_{{\Bbb R}^n} e^{zq(x)} e^{-\|x\|^2/2} \ dx = \det^{\quad -{1 \over 2}}\left(I - z Q\right),$$
where we choose the branch of  $\det^{ -{1 \over 2}}\left(I - z Q\right)$, which is equal to $1$ when $z=0$.
\endproclaim
\demo{Proof} The formula is well-known for real $z$. Since both sides of the identity are analytic functions of $z$, it holds for complex $z$ as well.
{\hfill \hfill \hfill} \qed
\enddemo

\subhead (3.3) Proof of Proposition 2.2 \endsubhead Using Lemmas 3.1 and 3.2, we write 
$$\split &{1 \over (2 \pi)^{n/2}} \int_{{\Bbb R}^n} e^{-\|x\|^2/2} \prod_{i=1}^m {\sin^2 \bigl(z q_i(x)\bigr) \over z^2 q_i^2(x)} \ dx \\&= 
{2^m \over (2 \pi)^{n/2}} \int_{{\Bbb R}^n} e^{-\|x\|^2/2} \left( \int_{[0,1]^m} \prod_{i=1}^m \left(1-t_i\right) \cos \left( 2z t_i q_i(x)\right) \ dt \right) \ dx \\
&={1 \over (2 \pi)^{n/2}} \int_{{\Bbb R}^n} e^{-\|x\|^2/2} \Biggl( \int_{[0,1]^m} \prod_{i=1}^m \left(1-t_i\right) \\ &\qquad \qquad \qquad \times \left( e^{2 z t_i \sqrt{-1} q_i(x)} + e^{-2z t_i \sqrt{-1} q_i(x)} \right)  \ dt \Biggr) \ dx  \\
&={1 \over (2 \pi)^{n/2}} \int_{{\Bbb R}^n} e^{-\|x\|^2/2} \Biggl( \int_{[0,1]^m} \prod_{i=1}^m \left(1-t_i\right) \\ &\qquad \qquad \qquad \times \sum_{\sigma_1, \ldots, \sigma_m=\pm 1} \exp\left\{ 2z \sum_{i=1}^m \sigma_i t_i 
\sqrt{-1} q_i(x)  \right\} \ dt \Biggr) \ dx \\
&=\int_{[0,1]^m} \Biggl( \prod_{i=1}^m \left(1-t_i\right) {1 \over (2 \pi)^{n/2}} \int_{{\Bbb R}^n} e^{-\|x\|^2/2} \\ 
 &\qquad \qquad \qquad \times  \sum_{\sigma_1, \ldots, \sigma_m=\pm 1} \exp\left\{ 2z \sum_{i=1}^m \sigma_i t_i \sqrt{-1} q_i(x) \right\} \ dx \Biggr) \ dt \\
 &=\int_{[0,1]^m} \prod_{i=1}^m \left(1-t_i\right) \sum_{\sigma_1, \ldots, \sigma_m=\pm 1} \det^{\qquad -{1\over 2}}\left(I -2z\sum_{i=1}^m  \sigma_i t_i \sqrt{-1} Q_i \right)  \ dt \\
 &=\int_{[-1,1]^m} \prod_{i=1}^m \left(1-|t_i| \right) \det^{\qquad -{1\over 2}}\left(I -2z\sum_{i=1}^m  t_i \sqrt{-1} Q_i \right)  \ dt
\endsplit$$
and the proof follows. 
{\hfill \hfill \hfill} \qed 

\head 4. Proof of Theorem 2.3 \endhead

The proof is very similar to that of Theorem 3.1 from \cite{BR19}.

It is more convenient to work with probability measures, as opposed to general measures. We start with a simple lemma, which provides a lower bound on the absolute value of the expectation of a complex-valued random variable. In what follows, we measure angles between non-zero complex numbers as between vectors in the plane, thus identifying ${\Bbb C}={\Bbb R}^2$.
\proclaim{(4.1) Lemma} Let $\Omega$ be a probability space and let $f: \Omega \longrightarrow {\Bbb C}$ be a random variable. Suppose that $f(\omega) \ne 0$ for all $\omega \in \Omega$, and, moreover, 
for any two $\omega_1, \omega_2 \in  \Omega$, the angle between $f(\omega_1) \ne 0$ and $f(\omega_2) \ne 0$ does not exceed $\theta$ for some $0 \leq \theta < 2\pi/3$. Then 
$$\left| \EE f \right| \ \geq \ \cos(\theta/2) \EE |f|.$$
\endproclaim
\demo{Proof} This is Lemma 3.3 from \cite{BR19}. 
{\hfill \hfill \hfill} \qed
\enddemo

\definition{(4.2) Definition} Let us fix some $0 \leq \theta_1, \ldots, \theta_m < 2\pi/3$. We denote by \newline $\GG(\theta_1, \ldots, \theta_m)$ the set of functions
$$g(t_1, \ldots, t_m) = \sum_{a \in {\Bbb Z}^m_+} c_a \ttt^a \tag4.2.1$$
where $c_a \in {\Bbb C}$ satisfy the inequalities
$$\sum\Sb a \in {\Bbb Z}^m: \\ i \in \supp a \endSb |c_a| \prod_{j \in \supp a} {1 \over \cos (\theta_j/2)} \ \leq \ {\theta_i \over 2} \quad \text{for} \quad i=1, \ldots, m.$$
\enddefinition

Some observation are in order. First, for every $g \in \GG(\theta_1, \ldots, \theta_m)$, the series (4.2.1) converges absolutely and uniformly on $[-1, 1]^m$ and we identify the series with a function 
$g: [-1, 1]^m \longrightarrow {\Bbb C} $.
Second, the set $\GG\left(\theta_1, \ldots, \theta_m\right)$ is convex. Third, let $I \subset \{1, \ldots, m\}$ be a set, let $g \in \GG(\theta_1, \ldots, \theta_m)$ and let $\widehat{g}$ be the function obtained from $g$ by fixing 
values of $t_i \in [-1, 1]$ for $i \in I$. Hence $\widehat{g}$ is a function in $t_i$ for $i \in \overline{I}=\{1, \ldots, m\} \setminus I$ and we have $\widehat{g} \in \GG(\theta_i: \ i \in \overline{I})$. We consider $\widehat{g}$ as a function 
$\widehat{g}: [-1, 1]^{m -|I|} \longrightarrow {\Bbb C}$.

\subhead (4.3) Proof of Theorem 2.3 \endsubhead Scaling the measures $\nu_1, \ldots, \nu_m$, if necessary, we assume that $\nu_1, \ldots, \nu_m$ and hence $\nu$ are probability measures on $[-1, 1]^m$.
Our goal is to show that if $g \in \GG(\theta_1, \ldots, \theta_m)$ then 
$$\left| \EE e^{g} \right| \ \geq \ \left(\prod_{i=1}^m \cos (\theta_i/2) \right) \EE \left| e^{g}\right|.$$
For a function $f: [-1, 1]^m \longrightarrow {\Bbb C}$ and a subset $I \subset \{1, \ldots, m\}$, we denote by $\EE_I f$ the conditional expectation of $f$, that is, the function of $t_i$ for $i \notin I$ obtained from $f$ by integrating over the variables $t_i$ with $i \in I$. We will apply this construction to $f=e^g$, where $g \in \GG(\theta_1, \ldots, \theta_m)$. If $I=\{i\}$ consists of a single element, we use a shorthand 
$\EE_i$ instead of $\EE_{\{i\}}$. By induction on $m$, we prove the following three statements. 
\bigskip
{\sl Statement} $1.m$. Let $g_0, g_1 \in \GG(\theta_1, \ldots, \theta_m)$ be two functions, which differ in a single coefficient $c_a$, equal to $c_{a,0}$ in $g_0$ and to $c_{a, 1}$ in $g_1$. Then the angle between 
$\EE e^{g_0} \ne 0$ and $\EE e^{g_1} \ne 0$ does not exceed 
$$| c_{a, 0} - c_{a, 1}| \prod_{i \in \supp a} {1 \over \cos (\theta_i/2)}.$$
\bigskip
{\sl Statement} $2.m$. Let $g \in \GG(\theta_1, \ldots, \theta_m)$, let $1 \leq i \leq m$ and let $I=\{1, \ldots, m\} \setminus \{i\}$. Let $h=\EE_I e^g$, so $h$ is a function of $t_i$. Then, for any $t_{i,1}, t_{i,2} \in [-1, 1]$, the angle between $h(t_{i,1}) \ne 0$ and $h(t_{i, 2}) \ne 0$ does not exceed $\theta_i$.
\bigskip
{\sl Statement} $3.m$. Let $g \in \GG(\theta_1, \ldots, \theta_m)$, let $I \subset \{1, \ldots, m\}$ and let \newline $\overline{I} =\{1, \ldots, m\} \setminus I$. Then 
$$|\EE e^g | \ \geq \ \left(\prod_{i \in I} \cos{ \theta_i\over 2} \right) \EE_I \left| \EE_{\overline{I}} e^g \right|. $$
In particular, if $I =\{1, \ldots, m\}$, we have 
$$\left| \EE e^g \right| \ \geq \ \left( \prod_{i=1}^m \cos{\theta_i \over 2}\right) \EE \left| e^g \right|  > 0.$$
\bigskip
We start with proving Statement $2.1$. In this case,
$$g(t)=\sum_{a \in {\Bbb Z}_+} c_a t^a$$
is a univariate function and 
$$\sum_{a  \in {\Bbb Z}_+} |c_a| \ \leq \ (\theta_1/2) \cos(\theta_1/2).$$
Furthermore, we have $h(t)=e^{g(t)}$. Clearly, $h(t) \ne 0$ for all $t$. Denoting by $\Im\thinspace w$ the imaginary part of a complex number $w$, we observe that for $-1 \leq t_1, t_2 \leq 1$, the angle between $h(t_1)$ and $h(t_2)$ does not exceed
$$\left| \Im\thinspace g(t_1) \right| + \left|\Im\thinspace g(t_2)\right| \ \leq \ 2 \sum_{a \in {\Bbb Z}_+} |c_a| \ \leq \ \theta_1 \cos(\theta_1/2) \ \leq \ \theta_1$$
and Statement $2.1$ follows.

We prove that Statement $2.m$ implies Statement $3.m$.
Let  $g \in \GG\left(\theta_1, \ldots, \theta_m\right)$. We proceed by induction on $|I|$. Suppose first, that $|I|=1$, so that $I=\{i\}$. Let $h=\EE_{\overline{I}} e^g$, so $h$ is a function of $t_i$. By Statement $2.m$, we have that the angle between
$h(t_1) \ne 0$ and $h(t_2) \ne 0$, for any two $t_1, t_2 \in [-1, 1]$ does not exceed $\theta_i$. Applying Lemma 4.1, we conclude that 
$$\left| \EE e^g \right| = |\EE h | \ \geq \ \cos (\theta_i/2) \EE |h| = \cos (\theta_i/2) \EE_i \left| \EE_{\overline{I}} e^g \right|.$$
If $|I| > 1$, let us choose $i \in I$ and let $J=I \setminus \{i\}$. Then 
$$\left| \EE e^g \right| = \left| \EE_i \EE_J \EE_{\overline{I}} e^g \right|.$$
Let $g_t$ be the function obtained from $g$ by fixing $t_i=t$, so $g_t \in \GG\left(\theta_j: \ j \ne i \right)$ for all $t \in [-1,1]$. By the induction hypothesis, for all $t \in [-1, 1]$, we have 
$$\left| \EE e^{g_t}\right| =\left| \EE_J \EE_{\overline{I}} e^{g_t} \right| \ \geq \ \left(\prod_{j \in J} \cos(\theta_j/2)\right) \EE_J \left| \EE_{\overline{I}} e^{g_t} \right|.$$
On the other hand, by Statement $2.m$, for any two values of $-1 \leq t_1, t_2 \leq 1$, the angle between $\EE e^{g_{t_1}} \ne 0$ and $\EE e^{g_{t_2}} \ne 0$ does not exceed $\theta_i$ and hence by Lemma 4.1, 
we have
$$\split \left| \EE e^g \right| \ \geq \ &\cos(\theta_i/2) \int_{-1}^1 \left| \EE e^{g_t} \right| \ d \nu_i(t) \ \geq \ \left(\prod_{i \in I} \cos(\theta_i/2)\right) \EE_i \EE_J \left| \EE_{\overline{I}} e^g \right| \\
=& \left(\prod_{i \in I} \cos(\theta_i/2)\right) \EE_I \left| \EE_{\overline{I}} e^g \right| \endsplit $$
and Statement $3.m$ follows. 

Next, we prove that Statement $3.m$ implies Statement $1.m$.
Given $g_0, g_1 \in \GG(\theta_1, \ldots, \theta_m)$ which differ in a single coefficient $c_a$, for $0 \leq \beta \leq 1$, let 
$$g_{\beta}= (1-\beta) g_0 + \beta g_1,$$
so $g_{\beta}=g_0$ for $\beta=0$ and $g_{\beta}=g_1$ for $\beta=1$. We have $g_{\beta} \in \GG(\theta_1, \ldots, \theta_m)$. In particular, $\EE e^{g_{\beta}} \ne 0$ by Statement $3.m$ and hence we can choose a branch of the function $\beta \longmapsto \ln  \EE e^{g_{\beta}}$ for $0 \leq \beta \leq 1$. We have 
$$\split \ln \EE e^{g_1} - \ln \EE e^{g_0}= &\int_0^1 {d \over d \beta} \ln \EE e^{g_{\beta}} \ d \beta = \int_0^1 {(d/d\beta) \EE e^{g_{\beta}} \over \EE e^{g_{\beta}}} \ d \beta \\=
 &(c_{a,1}-c_{a,0}) \int_0^1 {\EE \left(\ttt^a e^{g_{\beta}} \right) \over \EE e^{g_{\beta}}} \ d \beta. \endsplit$$
 Let $I=\supp a$. Since $|t_i| \leq 1$ for all $i$, we have
 $$\left| \EE \left(\ttt^a  e^{g_{\beta}}\right)\right|=\left|\EE_I \left(\EE_{\overline{I}} \left(\ttt^a e^{g_{\beta}} \right) \right)\right| = \left|\EE_I \left( \ttt^a \EE_{\overline{I}} e^{g_{\beta}}\right)\right| \ \leq \ 
 \EE_I \left| \EE_{\overline{I}} e^{g_{\beta}} \right| .$$
 On the other hand, by Statement $3.m$, 
 $$\left| \EE e^{g_{\beta}}\right| =\left|  \EE_I \left(\EE_{\overline{I}} e^{g_{\beta}} \right)\right| \ \geq \ \left(\prod_{i \in I} \cos{ \theta_i\over 2} \right) \EE_I \left| \EE_{\overline{I}} e^{g_{\beta}} \right|. $$
 Summarizing,
 $$\left|\ln \EE e^{g_1} - \ln \EE e^{g_0}\right| \ \leq \ \left| c_{a,1} - c_{a, 0}\right| \prod_{ i \in I} {1 \over \cos(\theta_i/2)},$$
 and Statement $1.m$ holds.
 
 Finally, we show that Statement $1.(m-1)$ and  Statement $3.(m-1)$ imply Statement $2.m.$
 Without loss of generality, we suppose that $i=m$.
 Let us choose \newline $g \in \GG(\theta_1, \ldots, \theta_m)$, so 
 $$g(t_1, \ldots, t_m )=\sum\Sb a \in {\Bbb Z}^m_+ :\\ a=\left(\alpha_1, \ldots, \alpha_m\right) \endSb  c_a t_1^{\alpha_1} \cdots t_m^{\alpha_m}.$$
  For $t \in [-1, 1]$, let $g_t$ be the function obtained by fixing $t_m=t$ in $g$. Hence $g_t \in \GG\left(\theta_1, \ldots, \theta_{m-1}\right)$ 
  and $h(t)=\EE e^{g_t}$. We can write 
   $$g_t(t_1, \ldots, t_{m-1}) =\sum_{b \in {\Bbb Z}^{m-1}_+} \left(\sum\Sb a \in {\Bbb Z}^m_+: \\ a=(b, \alpha_m) \endSb c_a t^{\alpha_m}\right) \ttt^b.$$
 By Statement $3.(m-1)$, we have $h(t) \ne 0$ for all $-1 \leq t \leq 1$. As we change from $t=t_{m,1}$ to $t=t_{m,2}$, the coefficients of $g_t$ change, and 
 applying Statement $1.(m-1)$ repeatedly, we conclude that the angle between $h(t_{m,1}) \ne 0$ and $h(t_{m,2}) \ne 0$ does not exceed 
$$\split \sum_{b \in {\Bbb Z}^{m-1}_+}  \sum\Sb a \in {\Bbb Z}^m_+: \\ a=(b, \alpha_m), \alpha_m > 0 \endSb 2|c_a| \prod_{i \in \supp b} {1 \over \cos(\theta_i/2)} \ \leq \ &2\sum\Sb a \in {\Bbb Z}^m_+: \\ m \in \supp a \endSb
|c_a| \prod_{i \in \supp a} {1 \over \cos(\theta_i/2)} \\  \leq \ &\theta_m \endsplit$$
by the definition of $\GG(\theta_1, \ldots, \theta_m)$. Hence Statement $2.m$ holds.

This completes the induction. From Statement 3.m, we have 
$$\left| \EE e^g \right| \ \geq \ \left( \prod_{i=1}^m \cos (\theta_i/2) \right) \EE \left| e^g \right|.$$
{\hfill \hfill \hfill} \qed

\head 5. Proof of Theorem 2.1 \endhead

 Let
$$q_i(x)={1 \over 2} \langle Q_i x, x \rangle \quad \text{for} \quad i=1, \ldots, m,$$
where $Q_i$ are $n \times n$ symmetric matrices. Since 
$$\sum_{i=1}^m |q_i(x)| \ \leq \ \gamma \|x\|^2,$$
we have 
$$\left\| \sum_{i=1}^m z_i Q_i\right\| \ \leq \ 2 \gamma \ < \ {1 \over 2} \tag5.1$$
as long as $|z_1|, \ldots, |z_m| \leq 1$. By Proposition 2.2, we have 
$$\phi(z)=\int_{[-1, 1]^m}  \det^{\qquad-{1 \over 2}} \left(I - 2z \sqrt{-1} \sum_{i=1}^m t_i Q_i \right) \prod_{i=1}^m \left(1 - |t_i|\right) \ dt.$$
Next, we expand
$$\split &\det^{\qquad-{1 \over 2}} \left(I - 2z \sqrt{-1} \sum_{i=1}^m t_i Q_i \right)=\exp\left\{ -{1 \over 2} \ln \det \left(I - 2z \sqrt{-1} \sum_{i=1}^m t_i Q_i \right) \right\} \\
&\qquad = \exp\left\{ -{1 \over 2} \tr \ln \left(I - 2z \sqrt{-1} \sum_{i=1}^m t_i Q_i \right) \right\}\\&\qquad=\exp\left\{ {1 \over 2} \sum_{s=1}^{\infty} {(2z\sqrt{-1})^s \over s} \tr \left(\sum_{i=1}^m t_i Q_i \right)^s \right\} 
=\exp\left\{ \sum_{a \in {\Bbb Z}^m_+} c_a \ttt^a \right\},  \endsplit $$
where for $a=\left(\alpha_1, \ldots, \alpha_m \right)$, we have 
$$c_a={1 \over 2s} (2z\sqrt{-1})^s \sum_{(i_1, \ldots, i_s)} \tr \left(Q_{i_1} \cdots Q_{i_s}\right), \quad \text{where} \quad s=\alpha_1 + \ldots + \alpha_m$$
and the sum is taken over all sequences $(i_1, \ldots, i_s)$  that contain index $i$ exactly $\alpha_i$ times. 

We are going to apply Theorem 2.3 with $\nu_i=(1-|t_i|) \ dt_i$ and 
$$\theta_1 = \ldots = \theta_m = {\pi \over 2},$$
for which we bound
$$\sum\Sb a \in {\Bbb Z}^m_+: \\ i \in \supp a \endSb |c_a| \prod_{j \in \supp a} {1 \over \cos(\theta_j/2)} \ = \ \sum\Sb a \in {\Bbb Z}^m_+: \\ i \in \supp a \endSb |c_a| 2^{|\supp a| \over 2}.$$
To that end, we note that since $\rk Q_j \leq r$ and 
$\| Q_j \| \leq 2\gamma/r$ for all $j$, we have 
$$\left| \tr \left(Q_{i_1} \cdots Q_{i_s}\right) \right| \ \leq \ r \left({2 \gamma \over r}\right)^s.$$
Moreover, since for each matrix $Q_i$ there are at most $r$ matrices $Q_j$ (possibly with $j=i$) such that $Q_i Q_j \ne 0$. Hence, once $Q_{i_1}$ is chosen, there are at most $r^{s-1}$ choices of $Q_{i_2}, \ldots, Q_{i_s}$
with the property that $Q_{i_1} \cdots Q_{i_s} \ne 0$. Using that the trace of the product matrices is invariant under cyclic permutations, 
we conclude that 
$$\split \sum\Sb a \in {\Bbb Z}^m_+: \\ i \in \supp a \endSb |c_a| 2^{|\supp a| \over 2} \ \leq \  &{r \over 2} \sum_{s=1}^{\infty} 2^s 2^{s/2} \left({2 \gamma \over r}\right)^s r^{s-1} ={1 \over 2} \sum_{s=1}^{\infty}  (4 \sqrt{2} \gamma)^s  ={2 \sqrt{2} \gamma \over 1-4 \sqrt{2} \gamma} \\ < \ &{\pi \over 4}.\endsplit$$
Applying Theorem 2.3, we conclude that $\phi(z) \ne 0$, and, moreover,
$$| \phi(z)| \ \geq \ 2^{-m/2} \int_{[-1, 1]^m} \left| \det^{\qquad -{1 \over 2}} \left(I - 2z \sqrt{-1} \sum_{i=1}^m t_i Q_i \right)\right| \prod_{i=1}^m \left(1-|t_i|\right) \ dt.$$
By (5.1), the eigenvalues of the matrix 
$$\sum_{i=1}^m t_iQ_i$$
do not exceed $2\gamma$ in the absolute value, and hence 
$$(1-4 \gamma)^{-n/2} \ \geq \ \left| \det^{\qquad -{1 \over 2}} \left(I - 2z \sqrt{-1} \sum_{i=1}^k t_i Q_i \right)\right| \ \geq \ (1+4 \gamma)^{-n/2}.$$
The proof now follows.
{\hfill \hfill \hfill} \qed 

\head 6. Proof of Theorem 1.3 \endhead 

The polynomial interpolation method hinges on the following simple lemma.
\proclaim{(6.1) Lemma} Let $g(z)$ be a univariate polynomial of some degree $N > 0$ such that $g(z) \ne 0$ for all $z \in {\Bbb C}$ satisfying $|z| < \beta$, for some $\beta > 1$. 
Let us choose a branch of $h(z)=\ln g(z)$ for $|z| < \beta$ and let 
$$T_k(z)=h(0)+\sum_{i=1}^k {h^{(i)}(0) \over i!} z^i$$
be the Taylor polynomial of $h(z)$ of degree $k$, computed at $z=0$. Then 
$$\left| h(1) - T_k(1) \right| \ \leq \ {N \over (k+1) \beta^k (\beta -1)}.$$
\endproclaim
\demo{Proof} This is Lemma 2.2.1 of \cite{Ba16}.
{\hfill \hfill \hfill} \qed
\enddemo
As follows from Lemma 6.1, to approximate $h(1)$ within error $0 < \epsilon < 1$, it suffices to choose $k = O(\ln N - \ln \epsilon)$, where the implied constant in the ``$O$" notation depends on $\beta$ only. 
Furthermore, as is discussed in Section 2.2.2 of \cite{Ba16}, one can compute the values of $h(0), h'(0), \ldots, h^{(k)}(0)$ from the values of $g(0), g'(0), \ldots, g^{(k)}(0)$ in polynomial time. Hence to approximate
the value of $g(1)=e^{h(1)}$ within relative error $0 < \epsilon <1$, it suffices to compute the values of $g(0), g'(0), \ldots, g^{(k)}(0)$ for $k = O(\ln N - \ln \epsilon)$.  As long as computing $g^{(i)}(0)$ can be accomplished in $N^{O(i)}$ time, we get an algorithm of quasi-polynomial $N^{O(\ln N -\ln \epsilon)}$ complexity of approximating $g(1)$ within relative error $\epsilon$.

To prove Theorem 1.3, we fix a $\beta_1 > 1$ and let $\gamma_1=\gamma/\beta$, where $\gamma$ is the constant from Theorem 2.1. Hence we can choose $\gamma_1=0.09$
and we assume that 
$$|q_i(x)| \ \leq \ {\gamma_1 \|x\|^2 \over r} \quad \text{for} \quad i=1, \ldots, m.$$
It follows from Theorem 2.1 that for $|z| \leq \beta_1$ the integral
$$\phi(z)={1 \over (2\pi)^{n/2}} \int_{{\Bbb R}^n} e^{-\|x\|^2/2} \prod_{i=1}^m {\sin^2 \left(z q_i(x)\right) \over z^2 q_i^2(x)} \ dx, $$
converges absolutely to an analytic function, and that 
we have 
$$(1-4 \gamma)^{-n/2} \ \geq \ |\phi(z)| \ \geq \ 2^{-m/2} (1+4 \gamma)^{-n/2} \quad \text{for all} \quad |z| \leq \beta_1. \tag6.2$$
We are going to apply Lemma 6.1 to the Taylor polynomial $g=\phi_N(z)$ of $\phi(z)$, computed at $z=0$, and of a sufficiently large degree $N$.
First, we bound the coefficients of the Taylor expansion of $\phi(z)$.
\proclaim{(6.3) Lemma} Let
$$\phi(z)=1+ \sum_{k=1}^{\infty} f_k z^k \quad \text{for} \quad |z| < \beta_1$$
be the Taylor expansion of $\phi(z)$. Then 
$$|f_k| \ \leq \ \beta_1^{-k} (1-4 \gamma)^{-n/2} \quad \text{for} \quad k \geq 1.$$
\endproclaim
\demo{Proof} Let 
$${\Bbb S}^1 =\bigl\{z \in {\Bbb C}:\ |z| =\rho \bigr\}$$ be the circle of radius $\rho$ and let $\nu$ be the Haar probability measure on ${\Bbb S}^1$. Then we have  
$$\int_{{\Bbb S}^1} z^m \overline{z}^s \ d \nu =\cases \rho^{2m} &\text{if\ } s =m \\ 0 & \text{otherwise.} \endcases$$
Choosing an arbitrary $0 < \rho < \beta_1$, we obtain
$$f_k = \rho^{-2k} \int_{{\Bbb S}^1} \phi(z) \overline{z}^k \ d \nu,$$
from which, using (6.2), we obtain
$$|f_k| \ \leq \ \rho^{-k} \max_{|z|=\rho} |\phi(z)| \ \leq \ \rho^{-k} (1-4 \gamma)^{-n/2}.$$
Taking the limit as $\rho \longrightarrow \beta_1$, we conclude the proof.
{\hfill \hfill \hfill} \qed 
\enddemo

Let us choose $1 < \beta < \beta_1$, for example $\beta=(1+\beta_1)/2$. In view of Lemma 6.3, for a given $0 < \epsilon < 1$, we can choose 
$$N=N(\epsilon, \beta, n, m)=(n+m)^{O(\ln (n+m)-\ln \epsilon)}$$ so that 
$$\sum_{k=N+1}^{\infty} |f_k| \beta^k \ \leq \ {\epsilon (1-4\gamma)^{n/2} \over 3 \cdot 2^{m/2}}.$$
Then we consider the polynomial
$$\phi_N(z)=1+ \sum_{k=1}^N f_k z^k.$$
From (6.2) we conclude that 
$$\phi_N(z) \ne 0 \quad \text{if} \quad |z| <  \beta $$
and
$$|\phi_N(1) - \phi(1)| \ \leq \ (\epsilon/3)|\phi(1)|.$$
Finally, we use Lemma 6.1 with $g(z)=\phi_N(z)$ to approximate $\phi_N(1)$ within relative error $\epsilon/3$. The resulting number will approximate $\phi(1)$ within relative error $\epsilon$.

It remains to show that we can compute 
$$\phi^{(k)}(0)=k! f_k$$
in $(m+n)^{O(k)}$ time.
For a fixed $x \in {\Bbb R}^n$, let
$$p_i(z)={\sin^2 \bigl(z q_i(x)\bigr) \over z^2 q_i^2(x)} \quad \text{for} \quad i=1, \ldots, m.$$
By Lemma 3.1, 
$$p_i(z)=2 \int_0^1 (1-t) \cos\bigl(2 zt q_i(x) \bigr)\ dt = {\sin^2 (z q_i(x)) \over z^2 q_i^2(x)}.$$
Hence, by symmetry
$$p_i^{2k-1}(0)=0,$$
and
$$p_i^{(2k)}(0)=2\int_0^1 (1-t) (-2t q_i(x))^{k} \ dt.$$
Consequently,
$$\phi^{(2k-1)}(0)=0$$ and 
$$\split \phi^{(2k)}(0)=&\sum\Sb k_1, \ldots, k_s: \\ k_1 + \ldots + k_s=k \endSb {2^s \over (2 \pi)^{n/2}} \\
&\times \int_{{\Bbb R}^n} e^{-\|x\|^2/2} \left( \int_{[0, 1]^s} \prod_{i=1}^s (1-t_i) (-2 t_i q_i(x))^{k_i}  \ dt  \right) \ dx. \endsplit$$
The sum contains $m^{O(k)}$ summands, and each summand splits into the product of the integrals 
$$\int_0^1 (1-t_i) (-2 t_i)^{k_i} \ d t_i=(-2)^{k_i} \left( {1 \over k_i+1} - {1 \over k_i+2} \right)$$
and the integral 
$${1 \over (2 \pi)^{n/2}} \int_{{\Bbb R}^n} \prod_{i=1}^s q_i^{k_i} (x) \ dx. \tag6.4$$
The last integral can be computed by expanding 
$$\prod_{i=1}^s q_i^{k_i}(x) = \sum_{a \in {\Bbb Z}^n_+}  b_a \xx^a  \quad \text{where} \quad \xx^a =x_1^{\alpha_1} \cdots x_n^{\alpha_n} \quad \text{for} \quad a=\left(\alpha_1, \ldots, \alpha_n \right)$$
and $x=\left(x_1, \ldots, x_n \right)$,
and using that 
$${1 \over \sqrt{2 \pi}} \int_{-\infty}^{+\infty} x^{\alpha} e^{-x^2/2} \ dx = \cases 0 &\text{if $\alpha$ is odd} \\ (\alpha-1)!! &\text{if $\alpha$ is even.} \endcases$$
 The complexity of the resulting algorithm is $(m+n)^{O(k)}$.
 
Another way to compute (6.4) is via extracting the coefficient of $t_1^{k_1} \cdots t_s^{k_s}$ in the Taylor series expansion at $t_1=\ldots =t_s=0$ of 
$${1 \over (2 \pi)^{n/2}} \int_{{\Bbb R}^n} e^{-\|x\|^2/2} \exp\left\{ \sum_{i=1}^s t_i q_i(x) \right\} \ dx = \det^{\qquad -{1\over 2}}\left(I - \sum_{i=1}^s t_i Q_i \right),$$
as described in \cite{Ba93}. The advantage of the latter approach that it has polynomial complexity as long as $s$ remains fixed in advance.
{\hfill \hfill \hfill} \qed

\enddocument
\end